\documentclass{article}

\title{Markov chains with exponential return times are finitary}
\author{Omer Angel
  \and Yinon Spinka}
\date{July 2019}

\AtEndDocument{
  \bigskip
  \small
  \textsc{Omer Angel, Yinon Spinka} \par
  \textsc{Department of Mathematics, University of British Columbia} \par
  \textit{Email:} \texttt{\{angel,yinon\}@math.ubc.ca} \par
}

\usepackage{amsmath,amssymb,amsfonts,amsthm}
\usepackage[colorlinks=true,citecolor=black,urlcolor=blue,pdfborder={0 0 0}]{hyperref}
\hypersetup{linkcolor=[rgb]{0,0,0.6}}
\usepackage{graphicx}
\usepackage{enumitem}
\usepackage{color}
\usepackage{tikz}
\usepackage[margin=1cm]{caption}
\usepackage[protrusion=true
]{microtype}

\usepackage[nameinlink]{cleveref}
  \crefname{theorem}{Theorem}{Theorems}
  \crefname{thm}{Theorem}{Theorems}
  \crefname{mainthm}{Theorem}{Theorems}
  \crefname{lemma}{Lemma}{Lemmas}
  \crefname{lem}{Lemma}{Lemmas}
  \crefname{remark}{Remark}{Remarks}
  \crefname{prop}{Proposition}{Propositions}
  \crefname{defn}{Definition}{Definitions}
  \crefname{corollary}{Corollary}{Corollaries}
  \crefname{section}{Section}{Sections}
  \crefname{figure}{Figure}{Figures}

\newtheorem{thm}{Theorem}

\newtheorem{lemma}[thm]{Lemma}

\newtheorem{prop}[thm]{Proposition}
\theoremstyle{definition}

\newtheorem{remark}[thm]{Remark}
\newtheorem*{remark*}{Remark}

\newtheorem{question}[thm]{Question}

\renewcommand{\P}{\mathbb P}
\newcommand{\Z}{\mathbb Z}
\newcommand{\E}{\mathbb E}

\newcommand{\N}{\mathbb N}
\newcommand{\eps}{\varepsilon}

\newcommand{\cF}{\mathcal{F}}

\renewcommand{\P}{\mathbb{P}}
\newcommand{\1}{\mathbf{1}}

\newcommand{\iid}{i.i.d.}

\begin{document}

\maketitle

\begin{abstract}
  Consider an ergodic Markov chain on a countable state space for which the return times have exponential tails.
  We show that the stationary version of any such chain is a finitary factor of an i.i.d.\ process.
  A key step is to show that any stationary renewal process whose jump distribution has exponential tails and is not supported on a proper subgroup of $\Z$ is a finitary factor of an i.i.d.\ process.
\end{abstract}

\section{Introduction}

A fundamental problem in ergodic theory is to understand which processes can be expressed in terms of which other processes as factors.
A particular case where this question is not resolved in general is which processes can be expressed as finitary factors of \iid\ (independent and identically distributed) processes.
In this note we give a simple proof that certain basic processes can be expressed as finitary factors of \iid\ processes.

A process $X=(X_n)_{n\in\Z}$ is called a \textbf{factor} of a process $Y=(Y_n)_{n\in\Z}$ if $X$ can be written in the form $X_n = \varphi(\dots,Y_{n-1},Y_n,Y_{n+1},\dots)$ for some measurable function $\varphi$.
The factor is \textbf{finitary} if $\varphi$ has the following property: there exists an almost surely finite stopping time $R$ with respect to the filtration $(\cF_r)_{r \ge 0}$ generated by $\{Y_n\}_{|n| \le r}$ such that $X_0$ is measurable with respect to $\cF_R$. In other words, the value of $X_0$ can be determined just by looking at the finitely many variables $(Y_n)_{|n|\le R}$.
Such a stopping time $R$ is called a \textbf{coding window}.

Recall that a Markov chain on a countable state space is called ergodic if it is irreducible, aperiodic and positive recurrent.
Note that we include aperiodicity in the definition.
Such a Markov chain is said to have \textbf{exponential return times} if for some (and hence every) state $s$, the time $T$ to return to $s$ for the Markov chain started at~$s$ satisfies $\P(T\ge n) \leq C e^{-cn}$ for some $C,c>0$ (which may depend on the state~$s$) and all $n$.

\begin{thm}\label{thm:mc-is-finitary}
  Let $M=(M_n)_{n \in \Z}$ be a stationary ergodic Markov chain on a countable state space with exponential return times.
  Then $M$ is a finitary factor of an \iid\ process with a coding window which has exponential tails.
\end{thm}

This theorem is a corollary of \cref{thm:renewal-is-finitary} below, concerned with representing renewal processes as finitary factors of \iid\ processes.
A discrete \textbf{renewal process} is a $\{0,1\}$-valued random process $X=(X_n)_{n \in \Z}$ such that the distances between consecutive 1's are independent and identically distributed. We call this latter distribution the \textbf{jump distribution}.
Renewal processes appear naturally in Markov chains:
If $(M_n)_{n \in \Z}$ is a stationary Markov chain and $s$ is some state, then $(\1_{\{M_n=s\}})_{n \in \Z}$ is a stationary renewal process.

Let $T$ be a random variable taking values in the positive integers.
We say that $T$ is \textbf{non-lattice} if its support is not contained in a proper subgroup of $\Z$, or equivalently, if $\gcd\{t:\P(T=t)>0\}=1$. We say that $T$ has \textbf{exponential tails} if there exist constants $C,c>0$ such that $\P(T \ge t) \le Ce^{-ct}$ for all $t \ge 1$.
For such processes, we prove the following.

\begin{thm}\label{thm:renewal-is-finitary}
  Let $X=(X_n)_{n \in \Z}$ be a stationary renewal process whose jump distribution is non-lattice and has exponential tails. Then $X$ is a finitary factor of an \iid\ process with a coding window which has exponential tails.
\end{thm}

We remark that the assumptions of the theorem on the jump distribution are necessary even if one is willing to drop the conclusion about the coding window.

\begin{prop}\label{lem:necessary}
  If a stationary renewal process is a finitary factor of an \iid\ process, then its jump distribution is non-lattice and has exponential tails.
\end{prop}



\subsection{Prior results.}

Expressing Markov chains as finitary factors of \iid\ processes is not new.
For Markov chains with finite state spaces, Akcoglu, del Junco and Rahe \cite{akcoglu1979finitary} proved that any such ergodic chain is a finitary factor of any other, provided the latter has strictly higher entropy.
Harvey, Holroyd, Peres, and Romik~\cite{harvey2006universal} proved that this can be done with a coding window which has exponential tails.
However, both of these results rely crucially on the finiteness of the state space.
We remark that if the source process is allowed to have sufficiently high entropy, then coupling from the past yields a very simple construction of a factor map with exponential tail for the coding window (still for Markov chains with finite state space).

When the Markov chain has countable state space, Rudolph~\cite{rudolph1982mixing} showed that, under the assumptions of \cref{thm:mc-is-finitary}, $M$ is finitarily isomorphic to an \iid\ process.
There are two main advantages to our approach:
First, our construction is short and explicit, whereas \cite{rudolph1982mixing} relies on earlier and more abstract work \cite{rudolph1981characterization}.
Secondly, we get an explicit, exponential bound on the coding window, whereas previous arguments do not provide any bound.

Our result does not yield an isomorphism, but instead yields exponential tails for the coding window, making the two results incomparable.
Moreover, our proof is based on a probabilistic argument which provides a more explicit construction of the finitary factor.

\subsection{Outline}

We break the proof of \cref{thm:renewal-is-finitary} into three parts from which the theorem immediately follows.

The first part is to show that a large class of renewal processes are finitary factors of i.i.d.\ processes. In the case when the jump distribution is bounded, one may use coupling from the past to obtain a simple proof of the fact that such a renewal process is a finitary factor of an \iid\ process (see \cref{rem:bounded-T}). A similar argument applies to a larger class of renewal processes and yields the following.

\begin{prop}\label{prop:renewal-finitary}
  Let $X$ be a stationary renewal process with jump distribution~$T$.
  Suppose that $T$ is unbounded and satisfies
  \[ \liminf_{n \to \infty} \P(T=n \mid T \ge n) > 0 .\]
  Then $X$ is a finitary factor of an i.i.d.\ process with a coding window which has exponential tails.
\end{prop}

The second part is to show that one may reduce the problem for a renewal process with jump distribution $T$ to a renewal process with a modified jump distribution. To define this modified jump distribution, let $T_1,T_2,\dots$ be independent copies of $T$.
Let $\mu \in (0,1)$ and let $N \sim \text{Geom}(\mu)$ be independent of $\{T_n\}_n$.
Here we use the convention that $N$ takes values in the positive integers and $\P(N=n)=\mu (1-\mu)^{n-1}$ for $n \ge 1$. Define
\[ T^*_\mu := T_1 + \cdots + T_N .\]

\begin{lemma}\label{lem:diluted-renewal-process}
Let $T$ a jump distribution and let $\mu \in (0,1)$.
Let $X$ and $X^*$ be stationary renewal processes with jump distributions $T$ and $T^*_\mu$, respectively.
If $X^*$ is a finitary factor of an \iid\ process, then so is $X$. Moreover, if $X^*$ has a coding window with exponential tails, then so does $X$.
\end{lemma}

The third and final part is to show that one may always reduce to the case of \cref{prop:renewal-finitary} by replacing a given jump distribution $T$ with some other jump distribution $T^*_\mu$ having very regular tails.

\begin{lemma}\label{lem:exp-geom}
Suppose that $T$ is non-lattice and has exponential tails.
Then for any sufficiently small $\mu>0$, there exist $c>0$ and $\kappa>\nu > 1$ such that
\[ \P(T^*_\mu=n) = c \nu^{-n} + O(\kappa^{-n}) \qquad\text{as }n \to \infty .\]
\end{lemma}

\paragraph{Acknowledgments.}
We thank Tom Hutchcroft for suggesting the problem to us, Thomas Budzinski for helpful discussions, and Jeff Steif for comments on an earlier draft and on the problem's history.
This work was supported in part by NSERC of Canada.

\section{Proofs}

Before going into the proofs, we make a simple observation. A basic consequence of renewal theory is that if $X$ is a stationary renewal process with jump distribution $T$, then $T':= \min \{ i\ge 1 : X_i=1\} - \max\{i \le 0 : X_i=1\}$, the size of the block of $X$ containing the origin, is a size-biased version of $T$. In particular, $T$ has exponential tails if and only if $T'$ does. For this reason, the distinction between the two will not be important for our bounds on the tails of the coding window.

\subsection{Proof of Theorem~\ref{thm:mc-is-finitary}}

We now show that \cref{thm:mc-is-finitary} is a corollary of \cref{thm:renewal-is-finitary}.

Let $M$ be an ergodic Markov chain on a countable state space, and let $s$ be some state such that the first return time to $s$ has exponential tails. Let $X$ be the process given by $X_n := \1_{\{M_n=s\}}$.
Then $X$ is a stationary renewal process whose jump distribution is non-lattice (due to the fact that $M$ is apreriodic) and has exponential tails. Thus, by Theorem~\ref{thm:renewal-is-finitary}, $X$ is a finitary factor of an \iid\ process with a coding window which has exponential tails.
Now, given $X$, $M$ consists of independent excursions of the Markov chain (from $s$ to $s$) of given lengths.
From this it is easy to see that $M$ is a finitary factor of an \iid\ process. Moreover, using that the coding window for $X$ has exponential tails and that the jump distribution of $X$ has exponential tails, it is easy to see that the coding window for $M$ also has exponential tails (this relies on the simple fact that the composition of finitary factors whose coding windows have exponential tails is again such a factor; see, e.g., \cite[Lemma~9]{harel2018finitary}).
\qed

\subsection{Proof of Proposition~\ref{lem:necessary}}

\cref{lem:necessary} is simple to prove and we include the brief proof for completeness.

Let $X$ be a stationary renewal process and suppose that it is a finitary factor of an i.i.d.\ process $Y$.
Then there exists an event of the form $\{(Y_{-m},\dots,Y_m) \in A\}$, for some positive integer $m$ and some measurable set $A$, which has positive probability and on which $X_0=1$ almost surely.

Let us first show that $T$ is non-lattice.
Note that for $n>2m$, the event $\{X_0=X_n=1\}$ has positive probability, since it is implied by the event $\{(Y_{-m},\dots,Y_m) \in A,~(Y_{n-m},\dots,Y_{n+m}) \in A\}$. On the other hand, if $T$ was supported in $k\Z$ for some $k \ge 2$, then the two events $\{X_0=1\}$ and $\{X_n=1\}$ could not simultaneously occur unless $n \in k\Z$. This shows that $T$ is non-lattice.
 
Let us now show that $T'$, and consequently $T$, has exponential tails.
This is the case, since the process $\bar X_n := \1_{\{Y_{n-m},\dots,Y_{n+m}\in A\}}$ is $2m$-dependent and $\bar X_n=1$ implies $X_n=1$. 
\qed

\subsection{Proof of Proposition~\ref{prop:renewal-finitary}}

Given the stationary renewal process $X$, we consider an auxiliary process $Z$ defined as one plus the distance to the nearest 1 to the left (possibly itself), i.e.,
\[ Z_i := 1 + \inf \{k\geq0 : X_{i-k}=1\}. \]
Note that $Z_i=1$ if and only if $X_i=1$.
It follows from the renewal property of $X$ that $Z$ is a stationary Markov chain on the state space $\N$, with transition probabilities given by
\begin{align*}
p(n,n+1) &= \P(T>n \mid T \ge n), & 
p(n,1) &= \P(T=n \mid T \ge n),
\end{align*}
and all other transitions have probability $0$.
We shall construct below the Markov chain $Z$ as a finitary factor of \iid\ variables.
Since the process $(\1_{\{Z_i=1\}})_i$ coincides with $X$, this will show that $X$ is also a finitary factor of an \iid\ process. 

Let $Y=(Y_i)_{i \in \Z}$ be independent uniform $[0,1]$ random variables.
In terms of the Markov chain representation, we follow the rule that a chain which is at state $n$ at time $i-1$ moves to state~1 at time $i$ if $Y_i \le p(n,1)$ and to state $n+1$ otherwise.
This allows us to run (on the same probability space) the Markov chain starting at any time $t\in\Z$ and any state $s\in\N$. We denote such a copy of the chain by $Z^{(t,s)}=(Z^{(t,s)}_i)_{i \ge t}$.
We claim that as the starting time $t$ tends to $-\infty$, uniformly in $s$, the chains $Z^{(t,s)}$ almost surely converge to a limit process $Z$.
The limit is clearly a copy of the Markov chain with transitions governed by the \iid\ process $Y$.
Moreover, we shall show that $Z$ is a finitary factor of $Y$.

\medskip

For simplicity, let us consider first the case where $T$ has full support.
In this case $p(n,1)>0$ for all $n \in \N$, and by assumption, $\liminf_{n \to \infty} p(n,1)>0$.
Thus, there is some $a>0$ such that $p(n,1)>a$ for all $n$.
By the transition rule for $Z$, it follows that if $Y_i\le a$ then $Z^{(t,s)}_i=1$ for all starting times $t<i$ and all starting states $s$, and, in particular, $Z_i=1$.
Therefore, to determine the value of $Z^{(t,s)}_i$, it suffices to observe values of $Y$ to the left of $i$, up to the first place $j\le i$ where $Y_j\le a$.
Having determined that $Z_j=1$, one can compute $Z$ at the times in the interval $[j,i]$ from the process $Y$.
Hence, $Z_i$ is determined by the variables $Y_j,Y_{j+1},\dots,Y_i$.
Since $(\1_{\{Y_j \le a\}})_{j \in \Z}$ is a Bernoulli process, it follows that $Z$ is indeed a finitary factor of $Y$ with exponential tail for the coding window.

\medskip

\begin{figure}
  \centering
  \includegraphics[width=.85\textwidth]{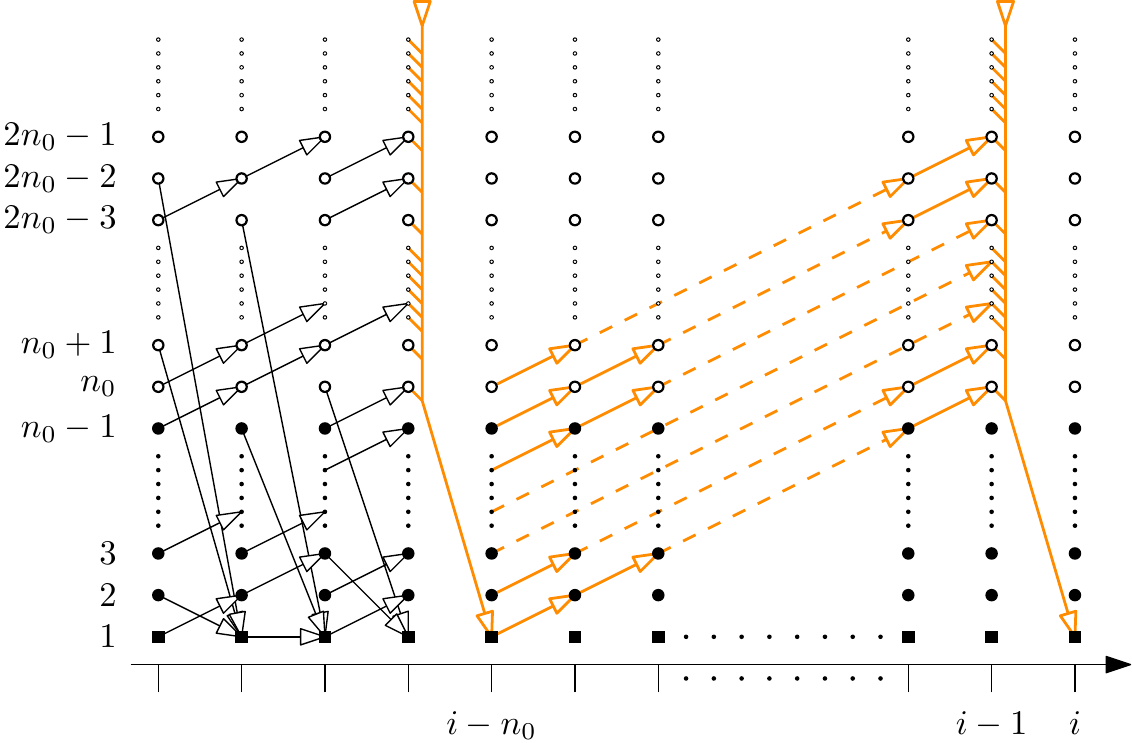}
  \caption{An illustration of the Markov chain and the event $E_i$. Some possible transitions between time $i-n_0-4$ and $i-n_0-1$ are shown. The occurrence of $E_i$ implies certain transitions between time $i-n_0-1$ and~$i$. These transitions are depicted in orange. Other transitions in this time frame are not depicted (and they are not relevant for the occurrence of $E_i$).}
 \label{fig:chain}
\end{figure}

It remains to deal with the case where $T$ does not have full support, i.e., when $p(n,1)=0$ for some $n$.
The main assumption of the proposition implies that there exists $n_0 \in \N$ such that
\[ a := \inf_{n \ge n_0} p(n,1) > 0 .\]
Since $T$ is unbounded, $p(n,n+1)>0$ for all $n \ge 1$. Thus,
\[ b := \max_{1 \le n \le 2n_0} p(n,1) < 1 .\]
For $i \in \Z$, consider the event
\[ E_i := \{ Y_{i-n_0} \le a \} \cap \{ Y_{i-n_0+1}, \dots, Y_{i-1} > b \} \cap \{ Y_i \le a \} .\]
We claim that if $E_i$ occurs for some $i$, then $Z^{(t,s)}_i=1$ for any $s$ and any $t<i-n_0$ (see Figure~\ref{fig:chain}).
Indeed, note that $Y_{i-n_0} \le a$ implies that if a chain is at some state $n \ge n_0$ at time $i-n_0-1$, then it is at state $1$ at time $i-n_0$. Of course, if a chain is at some state $n<n_0$ at time $i-n_0-1$, then it is at state $1$ or $n+1$ at time $i-n_0$. Thus, on the event $\{Y_{i-n_0} \le a\}$, all chains started before time $i-n_0$ are at some state in $\{1,\dots,n_0\}$ at time $i-n_0$. Next, note that $Y_{i-n_0+1} > b$ implies that if a chain is at some state $n \le n_0$ at time $i-n_0$, then it is at state $n+1$ at time $i-n_0+1$. Thus, on the event $\{Y_{i-n_0} \le a\} \cap \{ Y_{i-n_0+1} > b \}$, all chains started before time $i-n_0$ are at some state in $\{2,\dots,n_0+1\}$ at time $i-n_0+1$. In the same manner, we see that, on the event $E_i$, all chains started before time $i-n_0$ are at some state in $\{n_0,\dots,2n_0-1\}$ at time $i-1$.
Finally, since $E_i \subset \{Y_i \le a\}$, we see that, on the event $E_i$, all chains started before time $i-n_0$ are at state 1 at time $i$.

Hence, if $E_i$ occurs for some $i$, then we are guaranteed that $Z_i=1$.
Since $E_i$ almost surely occurs for infinitely many negative (and also positive) $i \in \Z$,
\[ I(j) := \max \{ i \le j : E_i\text{ occurs} \} \]
is almost surely finite for all $j \in \Z$.
Let $Z_j = Z^{(I(j),1)}_j$, i.e., the value at time $j$ of the chain started at time $I(j)$ at state 1.

We claim that this defines the chain as a finitary factor of $Y$ with exponential tail for the coding window.
Indeed, to determine $Z_j$, we only observe the values of $Y$ in the interval $[I(j)-n_0,j]$. 
Moreover, $I(j)$ is also determined by these same values of $Y$. Thus, $n_0-I(0)$ is a stopping time, and hence also a coding window for determining $X_0$.
This shows that $Z$ is a finitary factor of the i.i.d.\ process $Y$.
%
%
%
Finally, since $|I(0)| / (n_0+1)$ is stochastically dominated by a geometric random variable with parameter $\P(E_i) = a^2 (1-b)^{n_0-1}$, we conclude that the coding window has exponential tails.
\qed

\begin{remark}\label{rem:bounded-T}
When $T$ is bounded, a similar construction works, where instead of coupling the single-step transitions of chains at different states using a single variable $Y_i$ as we have done above, we simply let them evolve independently. Since there are finitely many possible states and the chain is ergodic, all states will eventually couple (with an exponential tail on the number of steps needed for this to happen). This is a special case of the method of coupling from the past~\cite{propp1996exact}.
\end{remark}

\subsection{Proof of Lemma~\ref{lem:diluted-renewal-process}}

Before proving the lemma, let us first note that that the opposite implication is simple -- namely, if $X$ is a finitary factor of an \iid\ process, then so is $X^*$.
This is so, because $X^*$ can be obtained as a dilution of $X$.
More precisely, if one changes each 1 in $X$ to a 0, independently, with probability $1-\mu$, then the resulting process has the same distribution as $X^*$.
Thus, $X^*$ is a finitary factor (in fact, a finitary factor with constant 0 coding window) of $(X,Y)$, where $Y$ is a Bernoulli \iid\ process independent of $X$.
In particular, if $X$ is a finitary factor of an \iid\ process, then so is $X^*$.

To establish the direction stated in the lemma, one essentially needs to reverse the above construction. In other words, given $X^*$, one needs to change some 0's into 1's in order to obtain a process with the same distribution as $X$.
The key observation, which can be seen from the above description, is that if one splits $\Z$ into blocks according to the 1's in $X^*$, then the restrictions of $X$ to different blocks are conditionally independent given $X^*$ (or equivalently, given the positions and lengths of the blocks).
Thus, in any such block, one may just sample $X$ given the block size, and do so independently for different blocks.
This can be done as a finitary factor, for example by using the value of $Y_s$ from an \iid\ process $Y$ (independent of $X^*$) to determine $X$ in a block of the form $[s,t]$.
If the coding window for $X^*$ has exponential tails, then, since the jump distribution of $X^*$ has exponential tails by \cref{lem:necessary}, we easily deduce that the coding window for $X$ also has exponential tails (this relies on the simple fact that the composition of finitary factors whose coding windows have exponential tail is again such a factor; see, e.g., \cite[Lemma~9]{harel2018finitary}).
\qed

\subsection{Proof of Lemma~\ref{lem:exp-geom}}

Consider the probability generating function of $T$,
\[ G(z) := \E[z^T] = \sum_{n=1}^\infty \P(T=n) z^n .\]
Since $T$ has exponential tails, this series has radius of convergence $a>1$.
Note that $G$ is strictly increasing and continuous on $[0,a)$, 
Thus, the limit $b := \lim_{z \to a^-} G(z)$ exists and satisfies $b \in (1,\infty]$, and $G^{-1} \colon [0,b) \to [0,a)$ is analytic, strictly increasing and satisfies $G^{-1}(1)=1$.
In particular, for any $\mu<1-\frac 1b$, there exists a $\nu \in (1,a)$ such that $(1-\mu)G(\nu) = 1$.
By the triangle inequality, for any $z$ we have $|G(z)| \le G(|z|)$.
Since $T$ is non-lattice, this is a strict inequality when $z\neq|z|$.
It follows that $(1-\mu) |G(z)| < 1$ whenever $|z| \le \nu$ and $z \neq \nu$. 

Consider also the probability generating function of $T^*_\mu$,
\[ F(z) := \E\big[z^{T^*_\mu}\big] = \E\big[ G(z)^N \big]
  = \frac{\mu G(z)}{1-(1-\mu) G(z)} .\]
It follows that for any $\mu<1-\frac1b$, the generating function $F$ is analytic on the disc of radius $\nu$ around $0$.
Fix such a $\mu$.
We claim that $F$ has a simple pole at $z=\nu$.
Indeed, this follows from the formula for $F$, since $G$ is analytic in a neighborhood of $\nu$, and $(1-\mu) G(\nu) = 1$ and $G'(\nu) = \sum_{n=1}^\infty n \P(T=n) \nu^{n-1} \neq 0$.

We claim that, for some $\eps>0$, the only singularity of $F$ in the disc of radius $\nu + \eps$ around 0 is the simple pole at $\nu$.
Indeed, since $G$ is analytic in the disc of radius $a$ around 0 and since $\nu<a$, the only singularities of $F$ in the disc of radius $a$ are at points where $G(z)=\frac{1}{1-\mu}$.
One such point is at $z=\nu$.
By analyticity, there is a small open neighborhood $U$ of $\nu$ in which $G$ does not equal $\frac 1{1-\mu}$ except at $\nu$.
Since $|G(z)|<\frac 1{1-\mu}$ on $\{ z : |z| \le \nu \} \setminus \{\nu\}$, it follows from the compactness of $\{ z : |z| \le \nu \} \setminus U$ that there is an open set $V$ containing $\{ z : |z| \le \nu \} \setminus U$ on which $|G(z)|<\frac 1{1-\mu}$. Thus, the only solution to $G(z)=\frac{1}{1-\mu}$ in $U \cup V$ is $z=\nu$. Since $U \cup V$ is an open set containing $\{z : |z| \le \nu\}$, the required $\eps$ exists.

Let $c$ be minus the residue of $F$ at the simple pole at $\nu$.
Then the function
\[ H(z) := F(z) - \frac {c}{\nu-z} \]
is analytic on a disc of radius $\nu + \eps$ around $0$.
In particular, $H$ has a Taylor expansion around~0, $H(z)=\sum_{n=0}^\infty c_n z^n$, with radius of convergence at least $\nu + \eps$. Hence, $\limsup_{n \to \infty} |c_n|^{1/n} \le \frac1{\nu + \eps}$. Since $\frac c{\nu-z}$ has a Taylor expansion, $\frac c\nu \sum_{n=0}^\infty \nu^{-n} z^n$, with radius of convergence $\nu$, we see that $F$ has a Talyor expansion, $F(z) = \sum_{n=0}^\infty (c \nu^{-n-1} + c_n) z^n$, with positive radius of convergence. In particular, since $F$ is the probability generating function of $T^*_\mu$,
\[ \P(T^*_\mu=n) = c\nu^{-n-1} + c_n = (\tfrac c\nu)\nu^{-n} + O(\kappa^{-n}) \qquad\text{as }n \to \infty ,\]
for any $\kappa<\nu+\eps$, as required. Note also that $c>0$, since $\P(T^*_\mu=n)$ is asymptotic to $c \nu^{-n-1}$.
\qed

\section{Open problems}

We believe that the conclusion of \cref{thm:mc-is-finitary} may be strengthened by also requiring that the \iid\ process have entropy at most $\eps$ larger than that of the Markov chain $M$. We raise the following question of whether this can be done with no gap in the entropy.

\begin{question}
Let $M=(M_n)_{n \in \Z}$ be a stationary ergodic Markov chain on a countable state space with exponential return times.
Does there exist an \iid\ process $Y$ with the same entropy as $X$ such that $M$ is a finitary factor of $Y$ with a coding window which has exponential tails?
\end{question}

We remark that even for Markov chains on a finite state space, this does not follow from existing results.

\bibliography{library}

\end{document}